\documentclass{amsart}
\usepackage{amsfonts}
\usepackage{amsmath}
\usepackage{graphics}
\usepackage{indentfirst,graphicx,epsfig}
\usepackage{cite}

\newcommand{\rem}[1]{}

\numberwithin{equation}{section}

\def\eqn{\begin{equation}}
\def\enn#1{\label{#1} \end{equation}}

\def\f#1{Figure {\ref{#1}}}
\def\t#1{Table {\ref{#1}}}

\begin{document}

\title[KdV-BURGERS  TYPE EQUATION WITH FAST DISPERSION AND SLOW
DIFFUSION] {ANALYSIS AND COMPUTATION OF A DISCRETE
 KdV-BURGERS  TYPE EQUATION WITH FAST DISPERSION AND SLOW
DIFFUSION }

\date{August 19, 2009}

\author[Z. Artstein]{Zvi Artstein}
\address[Z. Artstein]
{Department of Mathematics \\
The Weizmann Institute of Science\\ Rehovot 76100, Israel}
\email{zvi.artstein@weizmann.ac.il}

\author[C.W. Gear]{C. William Gear}
\address[C.W. Gear]
{ Department of Chemical Engineering\\
Princeton University\\
Princeton, New Jersey 08544, USA}
\email{wgear@princeton.edu}

\author[I.G. Kevrekidis]{Ioannis G. Kevrekidis}
\address[I.G. Kevrekidis]
{ Department of Chemical Engineering and PACM\\
Princeton University\\
Princeton, New Jersey 08544, USA}
\email{yannis@princeton.edu}

\author[M. Slemrod]{Marshall Slemrod}
\address[M. Slemrod]
{ Department of Mathematics\\
University of Wisconsin--Madison\\
Madison, Wisconsin 53706, USA}
\email{slemrod@math.wisc.edu}

\author[E.S. Titi]{Edriss S. Titi}
\address[E.S. Titi]
{ Department of Computer Science and Applied Mathematics \\
Weizmann Institute of Science  \\
Rehovot 76100, Israel. {\bf Also:} Department of Mathematics
 and Department of Mechanical and  Aerospace Engineering \\
University of California \\
Irvine, CA  92697-3875, USA.} \email{edriss.titi@weizmann.ac.il}

\begin{abstract}
The long time behavior of the dynamics of a fast-slow system of
ordinary differential equations is examined. The system is derived
from a spatial discretization of a Korteweg-de Vries-Burgers type
equation, with fast dispersion and slow diffusion. The
discretization is based on a model developed by Goodman and Lax,
that is composed of a fast system drifted by a slow forcing term. A
natural split to fast and slow state variables is, however, not
available. Our approach views the limit behavior as an invariant
measure of the fast motion drifted by the slow component, where the
known constants of motion of the fast system are employed as slowly
evolving observables; averaging equations for the latter lead to
computation of characteristic features of the motion. Such
computations are presented in the paper.
\end{abstract}

\maketitle

\vskip0.25in

\\


\section{Introduction}\label{S-1}

We consider the  slow-fast   system of ordinary differential
equations
\begin{equation}\label{Eq1.1}
\frac{dU_k}{ dt}+\frac{U_k(U_{k+1}-U_{k-1})}{2h} =
   \varepsilon \frac{U_{k+1}- 2U_k + U_{k-1}}{ h^2},
\end{equation}
$k=1,2,\cdots,N$, with periodic boundary conditions $U_k=U_{k+N}$,
and where $h > 0$ and $\varepsilon >0$ are small parameters, such
that $ \varepsilon \ll h$ in a sense that will be specified later.
It is an elementary observation that (\ref{Eq1.1}) is a spatial
discrete analogue of the classical Burgers' equation
\begin{equation}\label{Eq1.2}
u_t + uu_x= \varepsilon u_{xx},
\end{equation}
subject to periodic boundary conditions, with a basic domain
$[0,2\pi]$, and $h =\frac{2 \pi}{N}$ is the size of the uniform mesh
with $N$ points. But a more careful analysis (see, for example,
Goodman and Lax \cite{GoLa}) shows that a higher order approximation
of (\ref{Eq1.1}) is provided by the Korteweg-de Vries-Burgers' type
equation
\begin{equation}\label{Eq1.3}
u_t +  u (u_x+ \frac{h^2}{6} u_{xxx})= \varepsilon u_{xx},~~
 0\le x \le 2\pi,
\end{equation}
with periodic boundary conditions
\begin{equation}\label{Eq1.4}
u(0,t)=u(2\pi,t).
\end{equation}
In section 2 we establish a formal derivation of (\ref{Eq1.1}) from
(\ref{Eq1.3}), following Goodman and Lax \cite{GoLa}. This relation
of (\ref{Eq1.1}) to (\ref{Eq1.3}), under the assumption that $
\varepsilon \ll h$, reflects a competition between fast dispersion
and slow diffusion mechanisms. Our goal in this paper is to find an
effective equation governing the long time behavior of the ``slow''
diffusion, obtained via an ``averaging'' of the ``slow'' right hand
side of (\ref{Eq1.1}) against the fast oscillations originating from
the left hand side of (\ref{Eq1.1}). Intuition suggests that
averaging will be relevant; however, implementing classical
averaging results is not possible; we elaborate on this in the body
of the paper, but note here that the source of the difficulty is to
establish the existence and the identity of the ``angle variables"
in, say, a Hamiltonian representation, with respect to which
averaging should be carried out. Solutions of (\ref{Eq1.1}) can be
presented as images of solutions of a Hamiltonian system, but the
latter are not bounded, so do not induce relevant information (this
issue is discussed in section 3). Our approach relies on an
alternative new theory of computing slow observables for singularly
perturbed differential equations \cite{AKST}. This theory has the
distinct advantage of not requiring an explicit knowledge of the
``angle variables" and for this problem will produce a system of
equations for the slow observables of our approximate system of
ordinary differential equations. The Lax invariants, which are
traces of powers of a matrix component in a Lax pair, as developed
systematically in Goodman and Lax \cite{GoLa} and the related work
of Kac and van Moerbeke \cite{KavM}, are slow observables. We employ
these and other invariants of the fast flow in the computations of
the solutions of the discrete approximation.

The implementation of averaging techniques for systems of equations
which are separated into fast and slow components is well known.
These techniques in addition to component splitting often require
explicit information on the fast dynamics, e.g. time periodicity,
existence of limit cycles, or known stochastic behavior (see for
example the survey of Givon, Kupferman, and Stuart \cite{GKS}, see
also E et al. \cite{EEnLvE} and Kevrekidis et al. \cite{KGHKRT}).
The issue of averaging in the absence of explicit knowledge of the
fast dynamics was first considered by Artstein and Vigodner
\cite{AV} in 1996 and continued in Artstein \cite{A1}, \cite{A2},
Artstein and Slemrod \cite{AS1}, \cite{AS2}. A rigorous theory of
averaging in the absence of both a component split and explicit
limit characteristics of the fast dynamics was first given in
\cite{AKST}. The possible application of theory of \cite{AKST} to a
system with fast dispersion and slow diffusion was alluded to in
\cite{AKST} and this paper completes part of that program. In fact
when the theory is coupled with techniques such as projective
integration it provides, as we demonstrate in section 6, a fast and
accurate method for computation of fast dispersion systems
undergoing a slow perturbation. Of course our method does require
some special features of the fast system, i.e. a rich family of
first integrals of motion (invariants) of the type usually
associated with, but not limited to,  for instance, completely
integrable Hamiltonian systems. System (\ref{Eq1.1}) is particularly
interesting in that as far as available in the relevant literature,
e.g. Moser and Zehnder \cite{MZ}, there is no known one-to-one
correspondence between its fast part (when $\varepsilon=0$) and a
Hamiltonian system, nor do we know  whether solutions of this fast
system can be presented as  the image of solutions of a Hamiltonian
system with bounded dynamics. This fact renders the splitting into
fast and slow components difficult at best and perhaps impossible.

The paper has five sections after this Introduction. Section 2
derives the semi-discrete approximation to the system in a manner
suggested in the paper of Goodman and Lax \cite{GoLa}. It allows us
to represent, in Section 3, the fast part of the system in terms of
a Lax pair of matrices, thus producing a set of invariants of the
fast dynamics. We also comment there on the Hamiltonian-type
structure of the system. In section 4 we show that the positive
orthant is invariant under the dynamics of the full system
(\ref{Eq1.1}). Moreover, we use this fact to show the global
existence of solutions to (\ref{Eq1.1}) with initial data in the
positive orthant.  In Section 5 we apply the theory of \cite{AKST},
displaying Young measure solutions to the limit problem, and the
resulting evolution equations for the slow observables of our
equation.  Finally, in Section 6, we present numerical results for
the relevant ``effective slow'' dynamics obtained in the limit as
the parameter $\varepsilon \to 0$.


\vskip.1in


\section{Derivation of the semi-discrete system}  \label{S-2}

In this section we establish a formal derivation of (\ref{Eq1.1})
from (\ref{Eq1.3}) using the dispersive approximation scheme of
Goodman and Lax \cite{GoLa} suitably modified for our viscous
problem; for convenience of comparison we use the notation in
\cite{GoLa}.

For integers $k$ set $U_k(t)=u(kh,t)$. Hence periodicity of $u$ in
$x$ implies $U_k$ is periodic in $k$ with period $N$ where $N=
\frac{2\pi}{h}$. Assume $h$ is chosen so that $N$ is an even
integer. Taylor's theorem yields
\begin{equation}\label{Eq2.1}
{U_{k+1}-U_k} =
hu_x(kh,t)+\frac{1}{2}h^2 u_{xx}(kh,t)+
  \frac{1}{6}h^3u_{xxx}(kh,t) + O(h^4)
\end{equation}
and
\begin{equation}\label{Eq2.2}
U_k - U_{k-1}= hu_x(kh,t)-\frac{1}{2}h^2 u_{xx}(kh,t)+
   \frac{1}{6}h^3u_{xxx}(kh,t) + O(h^4).
\end{equation}
Addition of (\ref{Eq2.1}) and (\ref{Eq2.2}) yields
\begin{equation}\label{Eq2.3}
{U_{k+1}-U_{k-1}} = 2hu_x(kh,t)+\frac{2}{6}h^3 u_{xxx}(kh,t)
    +O(h^4),
\end{equation}
while subtraction of (\ref{Eq2.2}) from (\ref{Eq2.3}) gives
\begin{equation}\label{Eq2.4}
{U_{k+1}-2U_k+U_{k-1}} = h^2 u_{xx}(kh,t)+O(h^4).
\end{equation}
Notice that the right hand sides of the latter two equations
correspond to the expressions which constitute the continuous
equation (\ref{Eq2.1}). Thus the continuous system (\ref{Eq1.3}) is
formally $O(h^2)$ equivalent to the semi-discrete, i.e., continuous
in time only, discrete and $N$-periodic in $k$, system
\begin{equation}\label{Eq2.5}
\frac{dU_k}{dt}+\frac{U_k(U_{k+1}-U_{k-1})}{2h} =
   \varepsilon \frac{U_{k+1}- 2U_k + U_{k-1}}{h^2},
\end{equation}
with initial condition
\begin{equation}\label{Eq2.6}
U_k(0)=u_0(kh),~ k = 1,\ldots, N
\end{equation}
(which is our original system (\ref{Eq1.1})). The $N$-periodicity
means that in (2.5), and in equations later on, we use $U_0 = U_N$
and $U_{N+1} = U_1$. It is important to notice, though, that here
the step size of the approximation and the dispersion coefficient
are changed in a correlated way  through the equality $hN = 2 \pi$;
thus, a finer grid, namely, a larger but fixed $N$, corresponds to a
smaller dispersive coefficient; yet we always take $\varepsilon \ll
h. $

Goodman and Lax \cite{GoLa} examined the limit as $h \to 0$ of the
non-viscous version of (\ref{Eq1.3}) via the analogous semi-discrete
approximation. They proved, in particular, that during the time
period prior to the formation of a shock in the KdV type equation,
the solutions of the discrete approximation converge uniformly to
the solution of the continuous one; but after the time instant when
a shock forms in the Burgers equation, the discrete approximation
generates oscillations. We consider equations
(\ref{Eq2.5})-(\ref{Eq2.6}) as the semi-discrete approximation of
the KdV-Burgers equation (\ref{Eq1.3}); we wish to study the latter
for $h$ fixed and as $\varepsilon \to 0$, namely, on long time
intervals. To this end we rewrite (\ref{Eq2.6}) in the time scale
$\tau= \varepsilon t $, yielding
\begin{equation}\label{Eq2.7}
\frac{dU_k}{d\tau}+\frac{U_k(U_{k+1}-U_{k-1})}{2h \varepsilon} =
   \frac{U_{k+1}- 2U_k + U_{k-1}}{h^2}.
\end{equation}
 Denote by $U$ the vector $(U_1,\cdots,U_N)$; system
(\ref{Eq2.7}) can be rewritten as:
\begin{equation}\label{Eq2.8}
\frac{dU}{d\tau}=\frac{1}{\varepsilon} F(U) +G(U),
\end{equation}
where
$$
F(U) =   \frac{-1}{2 h}  \left( \begin{array}{c}
                  U_1(U_2-U_N)\\
                   U_2(U_3-U_1) \\
                   \cdot \\
                   \cdot\\
                   \cdot\\
                   U_N(U_1- U_{N-1}) \end{array}  \right)
 \qquad {\mbox {and}} \qquad
 G(U) =    \frac{1}{h^2}  \left(
\begin{array}{c}
                  U_2- 2U_1 + U_N \\
                   U_3- 2U_2 + U_1 \\
                   \cdot \\
                   \cdot\\
                   \cdot\\
                   U_1- 2U_N + U_{N-1} \end{array}  \right).
$$
 Now, in (\ref{Eq2.8}) (equivalently in (\ref{Eq2.7})) we
identify two contributions: The ``fast part", namely the vector
field  $\varepsilon^{-1} F(U)$ (for $h$ fixed), and the slow part,
namely, $G(U)$, which corresponds to the diffusion.
\vskip 15pt
%


\section{Analysis of the fast equation}  \label{S-3}

In this section we examine the ``fast part" of the vector field
determined by (\ref{Eq2.7}), i.e. we consider (\ref{Eq2.7}) after
replacing its right hand side by zero. We seek to find invariants of
the motion. These are not affected when multiplying the  equation by
a constant. Hence we write now the fast part of equation
(\ref{Eq2.7}) suppressing the coefficient $(2h\varepsilon)^{-1}$;
this amounts to a change of time scale,
$\sigma=(2h\varepsilon)^{-1}\tau$, without affecting other features
of the dynamics; namely, we examine the equation
\begin{equation}\label{Eq3.1}
\frac{dU_k}{d\sigma}+ U_k(U_{k+1}-U_{k-1}) = 0,
\end{equation}
where $U_k(0) > 0$ and $k = 1,2,\ldots,N$ with $N$ even and with
periodic boundary conditions (in particular, in (\ref{Eq3.1}) we
interpret $U_{N+1}= U_1$ and $U_0 = U_N$).
\vskip 6pt
{\bf Observation 3.1.}{\it The product $U_1 U_2\cdots U_N$  is an
invariant of equation (\ref{Eq3.1}); also the product  $U_2
U_4\cdots U_N$ of the even coordinates and the product  $U_1
U_3\cdots U_{N-1}$ of the odd coordinates of the state, are
invariants of equation (\ref{Eq3.1}). (We, actually, arrived at the
above quantities  guided by numerical simulations).}

{\bf Proof.} Straightforward, by showing that the time derivative of
each of the terms is equal to zero (in fact, establishing invariance
for two of the quantities implies that the third is also invariant).

\vskip 6pt
{\bf Corollary 3.2.} {\it The positive orthant, determined by the
condition $U_k > 0$ for all $k$, is invariant under the dynamics of
(3.1) (later in the section we establish a stronger property).}
\vskip 6pt
In the closing paragraph of this section we identify
Hamiltonian-type properties of the fast system (\ref{Eq3.1}). Now we
identify $\frac{N}{2}+1$ first integrals of the fast dynamics,
namely, $\frac{N}{2}+1$ constants of motion. To this end we review
the results of Goodman and Lax \cite{GoLa} for the fast system,
equation (\ref{Eq3.1}), and also borrow extensively from the
presentation of Kac and van Moerbeke \cite{KavM} and Moser
\cite{Mos1}, \cite{Mos2}; as far as the authors can tell the first
representation of (\ref{Eq3.1}) in terms of Lax pairs was given by
Moser \cite{Mos1}, but, as was noted by Moser himself, it was
analogous to that of Flaschka \cite{Fla}. Moser, however, imposes
Dirichlet boundary conditions while here we are interested in
periodic boundary conditions; this case was also discussed in Moser
and Zehnder \cite{MZ}.

We make a change of variables
\begin{equation}\label{Eq3.2}
A^2_k=U_k.
\end{equation}
With the additional change of time scale $s= -\frac{\sigma}{2}$, the
system (\ref{Eq3.1}) is transformed into
\begin{equation}\label{Eq3.3}
\frac{dA_k}{ds}=A_k(A_{k+1}^2-A_{k-1}^2),~ k=1,\ldots,N,
\end{equation}
with periodic boundary conditions, hence in (\ref{Eq3.3}) $A_1 =
A_{N+1}$ and $A_0 = A_N$. We note that any positive solution of the
system (\ref{Eq3.1}) can be rewritten as a solution of
(\ref{Eq3.3}); by Corollary 3.2 the product $A_1 A_2\ldots A_N$ is
also an invariant of the equation. Hence the positive orthant in
$R^N$ is invariant with respect to (\ref{Eq3.3}), and equations
(\ref{Eq3.1}) and (\ref{Eq3.3}) are equivalent for solutions with
positive initial data.

Now set $A=(A_1,\ldots,A_N)$, an $N$-dimensional vector of real
numbers, and extend $A_k$ periodically by $A_{k+N}=A_k$.

Define the operation $T$ on elements $A$ as the translation
$(TA)_k=A_{k+1}$ and use the abbreviation $TA=A_{+}$; likewise,
$T^{-1}$ is the translation in the opposite direction and
$T^{-1}A=A_{-}$.

Following Goodman and Lax \cite{GoLa}, Kac and van Moerbeke
\cite{KavM} and Moser \cite{Mos1}, \cite{Mos2}, we associate with
each vector $A$ the operators
\begin{equation}\label{Eq3.4}
L=AT+A_-T^{-1},
\end{equation}
where $A$ acts as a multiplication operator, i.e., when $T$ is
written in its matrix representation then $AT$ means that the $i$-th
row of $T$ is multiplied by the corresponding element, namely $A_i$,
of $A$. The matrix $L$ is symmetric with diagonal elements equal to
$0$. For the convenience of the reader and further reference we give
below the explicit terms when $N =6$; it should make the general
case transparent.

\begin{equation}\label{Eq3.5}
L = \left(\begin{array}{cccccc} 0&A_1&0&0&0&A_6 \\
               A_1&0&A_2&0&0&0 \\
                0&A_2&0&A_3&0&0 \\
                0&0&A_3&0&A_4&0 \\
                0&0&0&A_4&0&A_5 \\
                A_6&0&0&0&A_5&0   \end{array}\right)
\end{equation}
\vskip 3pt
\noindent With the aid of $L$ we now find conserved quantities for
the equation (\ref{Eq3.1}). To this end define
\begin{equation}\label{Eq3.6}
B=AA_+T^2-A_-A_{--}T^{-2}~,
\end{equation}
which is anti-symmetric.
Here $A_{--}=T^{-1}A_-=T^{-2}A$. Furthermore, the commutator
$[B,L]=BL-LB$ is given by
\begin{equation}\label{Eq3.7}
[B,L]=A(A^2_+-A^2_-)T-A_-(A-A_{--})T^{-1}.
\end{equation}
Note that (\ref{Eq3.3}) is equivalent to
\begin{equation}\label{Eq3.8}
\frac{dL}{ds} = [B,L].
\end{equation}
Then, of course,
\begin{equation}\label{Eq3.9}
\frac{dL^2}{ds} = L \frac{dL}{ds} + \frac{dL}{ds}L=[B,L^2]
\end{equation}
and more generally,
\begin{equation}\label{Eq3.10}
\frac{dL^p}{ds} =[B,L^p]
\end{equation}
for any natural number $p$. Since the right hand side of
(\ref{Eq3.10}) is a commutator its trace is equal to zero.
Therefore, for each $p$ we have
\begin{equation}\label{Eq3.11}
tr \frac{dL^p}{ds} \equiv 0
\end{equation}
where $tr$ signifies the trace. For $p$ odd the trace of $L^p$ is
zero, so it does not reveal any information concerning equation
(\ref{Eq3.3}). But for $p$ even, $p=2,4,\ldots, N$, the term $tr
L^p$ is a nontrivial polynomial of degree $p$. Thus, we have arrived
at $\frac{N}{2}$ non-trivial conserved quantities, namely, first
integrals, of equation (\ref{Eq3.3}). For example, the first two
non-trivial traces for $N \ge 6$ are
\begin{equation}\label{Eq3.12}
trL^2=2\sum^N_{k=1} A^2_k\,
\end{equation}
and
\begin{equation}\label{Eq3.13}
tr L^4=\sum^N_{k=1}(2A^2_kA^2_{k+1}+(A^2_k+A^2_{k-1})^2).
\end{equation}
The last trace for the case $N=6$ is
\begin{equation}\label{Eq3.14}
 tr L^6= \sum^6_{k=1}(A_k^2(A_{k-1}^2+A_k^2+A_{k+1}^2)^2 +
 A_{k-1}^2 A_k^2 A_{k+1}^2) + 12A_1A_2A_3A_4A_5A_6.
\end{equation}
\vskip 6pt
{\bf Remark 3.3.} {\it As stated in Observation 3.1 the product of
the $U_k$ is conserved by the equation (\ref{Eq3.1}); by
(\ref{Eq3.13}) each of the $U_k^2$ stays uniformly bounded. Hence,
for an initial condition with positive coordinates, the quantities
$U_k$ along a trajectory are bounded and bounded away from zero; in
particular, the solution exists for all time.}
\vskip 6pt
{\bf Remark 3.4.} {\it Invoking the relation (\ref{Eq3.2}), the
$\frac{N}{2}$ invariants for (\ref{Eq3.3}) give  $\frac{N}{2}$
invariants for (\ref{Eq3.1}). It is not easy to verify analytically
that the $\frac{N}{2}$ invariants are functionally independent,
namely, that they identify distinct $\frac{N}{2}$ dynamics;
numerical simulations show, indeed, that their gradients are
linearly independent, indicating their functional independence. In
addition, and as noted in Observation 3.1, the following  quantity
is also an invariant of (\ref{Eq3.1}),
\begin{equation}\label{Eq3.15}
U_1U_2\cdots U_{N}\,.
\end{equation}
It has been shown by simple algebra for $N = 6$ that the other
invariants mentioned in  Observation 3.1, namely

\begin{equation}\label{Eq3.16}
U_2U_4\cdots U_N \qquad {\rm and} \qquad  U_1U_3\cdots U_{N-1}
\end{equation}
are not independent of the $\frac{N}{2} + 1$ noted above and that
this is hypothesized to be true for all even $N$.  Altogether
 we get $\frac{N}{2} +1$   functionally independent invariants.
 We do not know whether there are any independent global (i.e. not
restricted to a specific lower dimensional manifold) first integrals
in addition to the $\frac{N}{2} +1$ that we have identified; as we
explain later on, the numerics indicates that the list is
exhausted.}

Considerations of the number of invariants may be significant as,
under an appropriate transformation, solutions of (\ref{Eq3.1}) are
differences of solutions of a completely integrable Hamiltonian
system. This was noted in Goodman and Lax \cite{GoLa}, referencing
the work of Kac and van Moerbeke \cite{KavM}, Flaschka \cite{Fla}
and Moser \cite{Mos1}, \cite{Mos2}, see also Moser and Zehnder
\cite{MZ}. The technique follows a change of variables
\begin{equation}\label{Eq3.17}
U_k = e^{z_k}
\end{equation}
under which the differential equation (\ref{Eq3.1}) takes  the form
\begin{equation}\label{Eq3.18}
\frac{dz_k}{d\sigma} = (e^{z_{k-1}} - e^{z_{k+1}}),
\end{equation}
which we relate  to the equation
\begin{equation}\label{Eq3.19}
\frac{dx_k}{d\sigma} = ( e^{x_{k-1} - x_k} + e^{x_k - x_{k+1}})
\end{equation}
via the relation
\begin{equation}\label{Eq3.20}
z_k = x_k - x_{k+1}.
\end{equation}
Observe that due to the periodicity of the $x_k$'s the sum of the
$z_k$'s is zero, which might be a source for an extra invariant.
Equation (\ref{Eq3.19}) is the well studied Toda lattice that forms
a completely integrable Hamiltonian system. It is clear, however,
from (\ref{Eq3.19}), that the solutions we are interested in,
namely, solutions of (\ref{Eq3.18}), are differences of unbounded
monotone increasing solutions (Moser \cite{Mos2} computes the
asymptotics of these solutions). Due to the unboundedness, this
Hamiltonian structure is not of much help in the computations of
solutions of (\ref{Eq3.1}). We do not know if, possibly, another
change of variables will induce a Hamiltonian structure on
(\ref{Eq3.1}) or will make solutions of (\ref{Eq3.1}) images of
bounded solutions of a Hamiltonian system.

%
\vskip 15pt
\section{A comment on the fast-slow system} \label{S-4}
\vskip 6pt
In this section we show that the positive orthant is invariant with
respect to the full fast-slow system (\ref{Eq1.1}).   Moreover, we
show that every such solution converges to a constant state
$U_1=U_2=\cdots=U_N$.

\vskip 6pt
{\bf Proposition 4.1.} {\it System (\ref{Eq1.1}) with $U_k(0) > 0$
for all $k=1,2,\ldots,N$, possesses a unique positive solution
$U_k(t) > 0$, for all $t \ge 0$. Moreover,
\begin{equation}\label{Eq4.1}
\prod_{k=1}^N U_k(t) \ge \prod_{k=1}^N U_k(0),
\end{equation}
and $U_k(t)$ converges to a constant, independent of $k$, for all
$k=1,2,\ldots,N$, as $t$ goes to infinity.}

{\bf Proof.} Since $U_k(0) > 0$, for all $k=1,2,\ldots,N$, then for
a short time we conclude from (\ref{Eq1.1}) that
\begin{eqnarray}
\frac{d}{dt} \sum_{k=1}^N \log U_k(t) &&= {\frac{\varepsilon}{
h^2}}\sum_{k=1}^N [\frac{U_{k+1}-U_k}{U_k}-\frac{U_{k}-U_{k-1}}{U_k}] \nonumber \\
 &&= \frac{\varepsilon}{h^2}\sum_{k=1}^N (U_{k+1}-U_k ) ( \frac{1}{U_k} - \frac{1}{U_{k+1}}) \label{Eq4.2} \\
&& \ge 0. \nonumber
\end{eqnarray}
The last inequality is true since the function $f(U)=\frac{1}{U}$ is
monotone decreasing. Hence,
\begin{equation}\label{Eq4.3}
\sum_{k=1}^N \log U_k(t) \ge \sum_{k=1}^N \log U_k(0),
\end{equation}
and
\begin{equation}\label{Eq4.4}
\prod_{k=1}^N U_k(t) \ge \prod_{k=1}^N U_k(0) > 0.
\end{equation}
Consequently, $U_k(t) >0$ for all $k=1,2,\ldots,N$, and for all
$t \ge 0$.

Since in addition
\begin{equation}\label{Eq4.5}
\sum_{k=1}^N U_k(t) = \sum_{k=1}^N U_k(0),
\end{equation}
we have $U_k(t)$ uniformly bounded and hence we have global unique
solutions. Finally, the expression $-\Sigma_{k=1}^N \log U_k(t)$ is
a Lyapunov function and via the well known LaSalle Invariance
Principle (see, for example, \cite{Ha}, Chapter 10, Theorem 1.3,
page 316) all solutions converge to the largest invariant set
contained in the set defined by the solution of the algebraic
equation
\begin{equation}\label{Eq4.6}
\sum_{k=1}^N (U_{k+1}-U_k ) ( \frac{1}{U_k} - \frac{1}{U_{k+1}})=0.
\end{equation}
Since the only possible solutions to the above equality are
$U_1=U_2=\cdots = U_N=const$, the proposition is proven.

%
\vskip 15pt
\section{ Young measures solutions and observables of the
fast-slow problem} \label{S-5}
\vskip 6pt
\noindent In this section we address the full system (\ref{Eq2.7}).
Let $T>0$ be arbitrary, but fixed. Based on the theory developed in
\cite{AKST} we represent limits, as $\varepsilon \to 0$, of
solutions as Young measures over the interval $[0,T]$; then examine
how slow observables of the system evolve over this arbitrary, but
fixed, interval of time. The analysis is the basis for the
computational developments in the closing section, section 6.

In the discussion that follows $N$ is fixed, hence $h$ is fixed. We
will write $U$ for the vector $(U_1,U_2,\ldots,U_k)$ in $R^N$. We
are interested then in limits of solutions $U_{\varepsilon}$ of
(\ref{Eq2.7}), as $\varepsilon \to 0$, over the interval $[0,T]$.

The Young measures we use are defined on a time interval $[0,T]$,
with values being probability measures on $R^N$, namely, on the
state space of solutions of (\ref{Eq2.7}). For the relevant theory
of Young measures see, e.g., \cite{AKST} and references therein. A
solution of (\ref{Eq2.7}), say $ U_{\varepsilon}(\tau) : [0,T] \to
R^N$, can be viewed as a Young measure when for each $\tau$ the
vector $ U_{\varepsilon}(\tau)$ is interpreted as a Dirac measure
supported on  $ \{U_{\varepsilon}(\tau)\}$. A Young measure can be
viewed also as a  measure on $[0,T] \times R^N$. Convergence of
Young measures of the full fast-slow system (\ref{Eq2.7}), as
${\varepsilon} \to 0$, is interpreted as the weak convergence of the
associated  measures on $[0,T] \times R^N$; the notion includes the
case where a sequence of point-maps into $R^N$, interpreted as Young
measures, converge to a Young measure. Convergence of solutions of
(\ref{Eq2.7}) in the sense of Young measures reflects that the
distribution of the location of the solutions converge. In
particular, the convergence in the sense of Young measures gives a
more complete description of the limit behavior of the full system
than the standard averaging method. This is because the Young
measures approach provides a mathematical framework, i.e. in the
sense of measures, for investigating and describing the limit
behavior of both the slow as well as fast dynamics, while the
averaging method restricts itself to the slow part of the dynamics.

The following result addresses the convergence of solutions in the
sense of Young measures; its proof is based on the developments in
previous sections and on  \cite{AKST}, Theorem 4.4.
\vskip 6pt
{\bf Theorem 5.1.} {\it Let the initial condition $U^0$ for
(\ref{Eq2.7}) be given (alternatively, let $U^0_\varepsilon$ be in a
compact set in the positive orthant of $R^N$). Let $T>0$ be given,
and fixed. Denote by $U_\varepsilon(\tau)$ the solution of
(\ref{Eq2.7}), for a given $\varepsilon$, over the interval $[0,T]$.
Then a subsequence $\varepsilon_k \to 0$ exists such that
$U_{\varepsilon_k}(\tau)$ converge as $k \to \infty$, in the sense
of Young measures on $[0,T]$, to a Young measure, say $\mu_0(\tau)$,
 $\tau \in [0,T]$. Moreover, for almost every $\tau \in [0,T]$
 the measure $\mu_0(\tau)$ is an invariant measure of the fast equation (\ref{Eq3.1}).}

{\bf Proof.} A straightforward computation shows that the quantity
$U_1 + \dots + U_N$ is not only an invariant of the fast flow
(\ref{Eq3.1}), see (\ref{Eq3.12}), but also of the full equation
(2.7) (indeed, it represents the total mass of the system, which,
under periodic boundary conditions, is preserved under viscosity).
Hence the family $U_\varepsilon(\tau)$ of functions, which is from
$[0,T]$ into the positive orthant of $R^N$, is bounded in $R^N$. In
particular, each solution in the sequence can be extended to the
entire interval $[0,T]$. Now the theorem follows from \cite{AKST},
Theorem 4.4.
\vskip 6pt
The next result addresses the evolution of observables of the
system. A measurement, or observable, is simply a function of $U$.
The notion of an orthogonal slow observable was defined in
\cite{AKST}; in our system it amounts to a measurement which is a
first integral, namely, an invariant, of the fast equation
(\ref{Eq3.1}). In Section 3 we identified $\frac{N}{2} + 1$ such
invariants. Denote by $v_j(U), j= 1,2,\ldots, \frac{N}{2}$ the
invariants which are the $\frac{N}{2}$ traces of $L^p,
p=2,4,\ldots,N$ (see Section 3, but note that in that section the
traces are expressed in terms of $A_k = U_k^{\frac{1}{2}}$); denote
by $v_{\frac{N}{2} + 1}(U)$ the invariant given in (\ref{Eq3.15}).
The invariance of each of the $v_j$ implies that whenever $\mu$ is
an invariant measure of (\ref{Eq3.1}), the observed value $v_j(U)$
is constant on the support of $\mu$. We denote this value by $\hat
v_j(\mu)$.
\vskip 6pt
{\bf Theorem 5.2.} {\it Let $U_{\varepsilon_k}(\tau)$ be as in the
statement of Theorem 5.1 which converge, as $k \to \infty$, in the
Young measures sense, to the Young measure $\mu_0(\tau)$, for $\tau
\in [0,T]$. Denote $\hat v_j(\tau)$ = $\hat v_j(\mu_0(\tau))$,
namely, the measurement on the limit invariant measures. Then for
each $v_j, j= 1,2,\ldots, \frac{N}{2}+1$, the function $\hat
v_j(\mu_0(\tau))$ satisfies the differential equation
\begin{equation}\label{Eq5.1}
\frac{d \hat v_j}{d\tau}(\tau) = \int_{R^N}
  (\nabla  v_j)(U)\cdot G(U)~\mu_0(\tau)(dU),~~~~\hat v_j(0) = v_j(U(0)),
\end{equation}
where $G(U) = G(U_1,\ldots,U_N)$ is the vector field given in
(\ref{Eq2.8}), and the $\nabla$ operator is with respect to the
vector $U$.
 Furthermore, the sequence of measurements
$v_j(U_{\varepsilon_k}(\tau))$ converge weakly to $\hat
v_j(\mu_0(\tau))$}, as $k \to \infty$, for $j= 1,2,\ldots,
\frac{N}{2}+1$.

{\bf Proof.} The claims form a particular case of Theorem 6.5 in
\cite{AKST}.
\vskip 6pt
Under stronger conditions on continuity properties of the limit
measure $\mu_0(\tau)$, the weak convergence in the previous result
becomes a strong convergence; we do not pursue this here since in
our system it is not easy to establish such continuity.

The differential equation (\ref{Eq5.1}) is not an autonomous
ordinary differential equation since the limit measure $\mu_0(\tau)$
is not determined by the observation $\hat v_j(\tau)$. In some cases
a set of observables determines the invariant measure; then, under
continuity assumptions one can write an ordinary differential
equation for this set of observables; without continuity a
differential inclusion would determine the limit evolution of the
observables. See  Theorem 6.9 and Remark 6.11 in \cite{AKST}. These
results can be stated in our case employing the terms in
(\ref{Eq5.1}), (\ref{Eq2.8}); we do not follow this direction since
we do not know to what extent the $\frac{N}{2} +1$ observables we
found form such a determining vector of observables. At any rate,
the structure described here allows one to make progress toward a
computational method, as displayed in the next section.
\vskip 15pt
%

\section{Computational study} \label{S-6}

In section we present  a computational study of  the theoretical
tools that were developed in \cite{AKST} and which were presented in
the previous sections for system (\ref{Eq2.7}) (or equivalently
(\ref{Eq2.8})), with small values of $\varepsilon$. Specifically, we
 introduce two efficient numerical schemes for computing the
evolution of the slow observables of (\ref{Eq2.7}), for small values
of $\varepsilon$. In particular, we   implement these schemes for
computing the evolution of the slow observables $v_j(U)$, $j=
1,2,\ldots, \frac{N}{2}+1$, that are discussed in section 5.
Moreover, we  demonstrate the computational efficiency of these
algorithms in comparison to direct numerical simulations of the full
system (\ref{Eq2.7}) for small values of $\varepsilon$. Both schemes
are based on numerical evaluation of the time derivatives of the
slow observables. That is, they are based on evaluating the time
derivatives of the invariants of the fast system (\ref{Eq3.1}),
${v}_j (U)$, $j=1,2,\cdots,\frac{N}{2}+1$, as they are drifted by
the slow field $G(U)$ in (\ref{Eq2.8}).

In addition to the evolution of the slow observables, the method can
simulate the evolution of the invariant measures of (\ref{Eq3.1})
themselves, thus getting a picture of the entire rapidly oscillating
trajectories of (\ref{Eq2.7}).

The first approach, the Young measures approach, is based on
computing the right hand side of (\ref{Eq5.1}). To this end, one has
to compute first the invariant measures $\mu_0(\tau)$ (see Theorem
5.1 and Theorem 5.2) of the fast system (\ref{Eq3.1}) by, for
example, integrating system (\ref{Eq3.1}) for an appropriate
interval of time, and then to advance ${\hat{v}}_j(\tau)$,  $j=
1,2,\ldots, \frac{N}{2}+1$, numerically, according to the numerical
evaluation of the right hand side of (\ref{Eq5.1}).

The second approach, the equation-free approach of \cite{KGHKRT} (it
is equation-free since it avoids the direct use of any analytic form
 for the equation of evolution of the slow observables $v_j$),
 is based on numerically evaluating the time derivative of
${v}_j(U_\varepsilon(\tau))$,  $j= 1,2,\ldots, \frac{N}{2}+1$, from
the detailed simulator of the solution $U_\varepsilon(\tau)$ of the
full fast-slow system (\ref{Eq2.7}), for small values of
$\varepsilon$. It is worth mentioning, however, that the success of
this method relies on the {\it a priori} knowledge that
${v}_j(U_\varepsilon(\tau))$, are slow observables of (\ref{Eq2.7})
and that at the limit, when $\varepsilon \to 0$, their evolution is
governed by a deterministic  equation; that may not be at hand.
These facts were established in the previous sections.

The efficiency of both algorithms stems from the fact one can take
relatively large time steps for advancing the slow observables,
${\hat{v}}_j(\tau)$ in the first approach, and
${v}_j(U_\varepsilon(\tau))$, for small values of $\varepsilon$, in
the second approach. We call this step in both algorithms  the
projective step.

It is worth mentioning the philosophical difference between the two
approaches. The first approach, which is based on Theorem 5.1 and
Theorem 5.2, takes advantage of the fact that the evolution of the
slow observables ${v}_j(U_\varepsilon(\tau))$, $j= 1,2,\ldots,
\frac{N}{2}+1$, along the solution $U_\varepsilon(\tau)$ of system
(\ref{Eq2.7}), are approximated, for small values of $\varepsilon$,
by their corresponding limits ${\hat{v}}_j(\tau)$, $j= 1,2,\ldots,
\frac{N}{2}+1$, as $\varepsilon \to 0$. Therefore, the first
approach focuses on computing the slow observables $\hat{v}_j(\tau)$
of the limit equation, as $\varepsilon \to 0$. This is because the
slow observables $\hat{v}_j(\tau)$ will then be  considered as the
approximations to the corresponding ``real" slow observables
$v_j(U_\varepsilon(\tau))$, for small values of $\varepsilon$. The
second approach, which is based on \cite{KGHKRT}, deals, on the
other hand, directly with system (\ref{Eq2.7}) without resorting to
the limit equation, equation (\ref{Eq5.1}), but it implicitly
assumes the existence of a limit equation.

Recall that in section 3 we  established that there are at least
$\frac{N}{2}+1$ invariants of the fast system (\ref{Eq3.1}),
therefore each trajectory is constrained to a manifold of dimension
at most $\frac{N}{2}-1$. This appears to be achieved for the case of
$N = 6$ as shown in \f{f1},
\begin{figure}[h]
\centerline{\psfig{figure=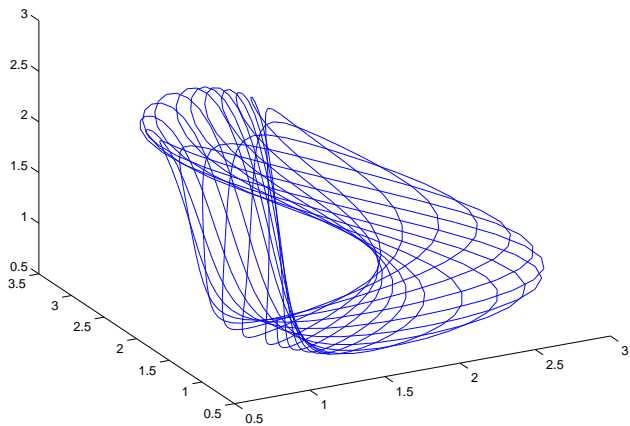,height=2.5in,height=2.5in}}
\caption{Torus for the case $N = 6$ of system (\ref{Eq3.1}). Initial
values were [1 1 1 3 2 1].}\label{f1}
\end{figure}
which illustrates the solution of (\ref{Eq3.1}) for one initial
value plotted in the first three variables, $U_1(\sigma),
U_2(\sigma)$ and $U_3(\sigma)$.  The solution lies on the surface of
a 2-torus in this 3D projection. The same observation holds for all
3D projections of this case except those projections into either of
all the odd-numbered or all the even-numbered coordinates (where it
necessarily lies on the surface of one sheet of a two sheeted
hyperboloid; this is because of the invariants noted in Observation
3.1). The trajectory ``orbits" around the ``long" direction of the
torus, precessing from orbit to orbit.
We assume, however, that we can approximate it, and the relevant
invariant Young measure, $\mu_0$,  by integrating over an
appropriate finite interval of time. Indeed, we have found that
integrating for a few ``orbits" over the ``long" direction of the
torus can give an adequate approximation for our application.

For $N = 6$ a typical solution of (\ref{Eq3.1}) was found to be
approximately a combination of two periodic motions, the fastest
corresponding to an orbit around the long direction of the torus and
the slow one corresponding to the precession time around the short
direction. See, for instance,  the graph of the first component,
$U_1(\sigma)$,  of the solution of system (\ref{Eq3.1}) shown in
\f{f1a}.
\begin{figure}[h]
\centerline{\psfig{figure=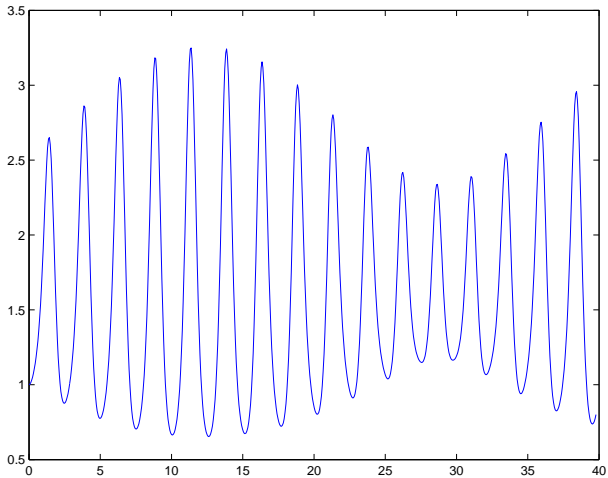,height=2.5in,height=2.5in}}
\caption{$U_1(\sigma)$ for the case in \f{f1}.}\label{f1a}
\end{figure}

Whereas the general theory, as described in section 3, shows that
there is a set of at least $\frac{N}{2}+1$ invariants for {\em all}
solutions of the $N$-dimensional problem (\ref{Eq3.1}), some
particular solutions may satisfy additional invariants which we
might call {\em local invariants}. That there are initial conditions
leading to local invariants is easy to see. For instance,  the case
$N = 3$ has 2 invariants ($U_1 + U_2 + U_3$ and $U_1U_2U_3$). The
case $N = 6$, the solutions  of (\ref{Eq3.1}) that have started with
$U_1 = U_4$, $U_2 = U_5$, and $U_3 = U_6$, will continue to maintain
those identities along the trajectory and will have the same
solution as the $N = 3$ problem with those starting conditions.
Hence its solution will lie on a one-dimensional manifold.  This is
shown in \f{f2}. (The invariants in this case are the 2 invariants,
for the $N = 3$ case, and the 3 equalities, $U_i = U_{i+3}$.) The
case $N = 12$, for example, that is known to have at least 7
invariants, will have certain initial conditions with 10 invariants
(corresponding to 2 copies of the $N = 6$ case) or 11 invariants
(corresponding to 3 copies of the $N = 4$ case, or 4 copies of the
$N = 3$ case).
\begin{figure}[h]
\centerline{\psfig{figure=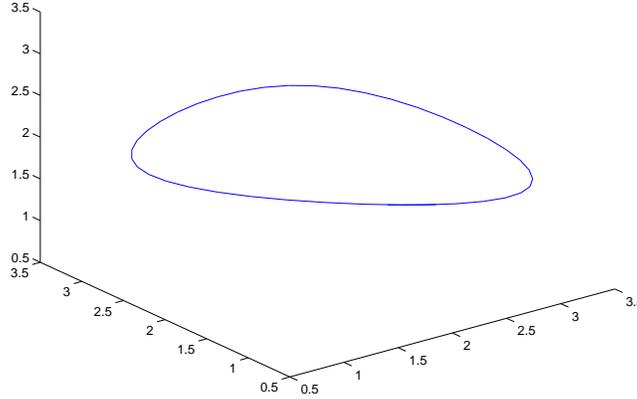,height=2.5in,height=2.5in}}
\caption{Case $N = 6$ as pair of $N = 3$ cases for system
(\ref{Eq3.1}). Initial values were [3 2 1 3 2 1].}\label{f2}
\end{figure}

\subsection{Two approaches for computing slow observables}
In this subsection we illustrate in more detail  the two approaches
for computing the behavior of the fast-slow system (\ref{Eq2.7})
which we rewrite in the form

\eqn \frac{dU_k}{dt} =  U_k(U_{k-1} - U_{k+1}) + \nu (U_{k+1} - 2U_k
+ U_{k-1}), \quad {\rm{for}} \quad k=1,2,\cdots, N, \enn{nuequ}
where $\nu= \frac{\varepsilon}{2h}$ is very small.

We do this for the case $N = 6$, but the technique is not restricted
to that case. Both approaches are based on evaluating the time
derivatives of the {\em slow observables}, $v_j(U)$,
$j=1,2,\cdots,\frac{N}{2}+1$, as we have discussed above.

In the first approach, which is based on the Young measures as
described in section 5, we utilize the right hand side of equation
(\ref{Eq5.1}) to evaluate the ``slow" derivatives of the slow
observables $\hat{v}_j (\tau)$, $j=1,2,\cdots,\frac{N}{2}+1$.  This
can be done by integrating the fast system (\ref{Eq3.1}), namely,
the $\nu = 0$ case of equation (\ref{nuequ}), to approximate the
invariant Young measure $\mu_0$ of (\ref{Eq3.1});  what we
practically  do to evaluate the right hand side of (\ref{Eq5.1}) is
to integrate over equally spaced time steps and approximate the
integral in the right hand side of (\ref{Eq5.1}) by an averaged
Riemann sum based on those points.

In the second approach, the equation-free method \cite{KGHKRT}, we
integrate the actual full fast-slow system equation (\ref{nuequ})
with small value of $\nu >0$, compute the slow observables
${v}_j(U_\nu(t))$,  $j=1,2,\cdots,\frac{N}{2}+1$, along the
trajectory $U_\nu(t)$, which are guaranteed to be slow by our
theory, and evaluate their slow derivatives by numerical
differencing.

In both approaches we use the time derivative evaluation of the slow
observables $v_j$, $j=1,2,\cdots,\frac{N}{2}+1$, in order to advance
them in time. We call this very last step  in our scheme the
``projective level". We observe that the gain in the speedup of both
schemes is due to the use of relatively large time steps at the
projective level of the schemes for advancing the slow obervables.
In order to repeat the process we need to initialize the detailed
simulator (i.e., find initial conditions, $\bar{U}$, for the
variables $U$ in the original $N$-dimensional space) consistent with
the ``projected" values of the $\frac{N}{2}+1$ slow observables.
That is, from the computed values of the slow observables, say
$\bar{v}_j$, $j=1,2,\cdots,\frac{N}{2}+1$, we need to find a
corresponding state variable $\bar{U}\in R^N$, for which
$v_j(\bar{U})=\bar{v}_j$, for all $j=1,2,\cdots,\frac{N}{2}+1$. This
procedure is called ``lifting" in the equation-free literature
\cite{KGHKRT}. In general this will involve the solution of a set of
$\frac{N}{2}+1$ nonlinear algebraic equations with $N$ unknowns (and
this implies that $\frac{N}{2}-1$ features of the solution can, in
principle, be selected (almost) arbitrarily; in our six-dimensional
example this ``arbitrary" selection was implemented by carrying the
values of two of the entries of the sought vector $\bar{U}$  over
from the last simulation.) As mentioned above, we take the point
$\bar{U}\in R^N$ as an initial value for the corresponding system,
i.e., an initial value for system (\ref{Eq3.1}) to compute the
corresponding Young measure $\mu_0$ in the first approach, or an
initial value for (\ref{nuequ}) to compute the corresponding
solution $U_\nu(t)$ in the second approach, and we repeat the above
described process. Note that if $\tilde{U} \in R^N$ is another point
that   has the same projected values of the slow observables as
those of $\bar{U}$, i.e. $v_j(\tilde{U})=v_j(\bar{U})= \bar{v}_j$,
for all $j=1,2,\cdots,\frac{N}{2}+1$, then we are not guaranteed
that the dynamics of the fast system (\ref{Eq3.1}) starting from
$\tilde{U}$ will  lead to the same invariant Young measure, $\mu_0$,
which is corresponding the trajectory starting from $\bar{U}$ (see
the relevant discussion regarding this matter in the closing
paragraph of section 5). This is unless the invariant Young measure
of (\ref{Eq3.1}) is unique, which is the case if the trajectory
corresponding to the solution of fast dynamics (\ref{Eq3.1})
``fills" the manifold defined by the invariants of (\ref{Eq3.1}). As
noted earlier, this manifold can be up to $\frac{N}{2}-1$
dimensional, since it is a sub-manifold of the manifold defined by
the intersection of the invariant manifolds $v_j(U)=constant$,
$j=1,2,\cdots,\frac{N}{2}+1$; although in reality it  might be of
lower dimension because of the additional {\em local invariants}.

\subsection{The first approach - using the invariant Young measures}

The invariants  $v_j(U)$, $j=1,2,\cdots,\frac{N}{2}+1$, of the fast
system (\ref{Eq3.1}) were identified in section 3 to be slow
observables of the full fast-slow system (\ref{Eq2.7}), or
equivalently of (\ref{Eq2.8}). Since these quantities  are constants
of motion for the fast dynamics, the contribution of the fast vector
field, $\varepsilon^{-1} F(U)$ in (\ref{Eq2.8}), to their evolution
is zero. Consequently, we  focus only on computing the evolution of
these slow observables as they are drifted by the slow vector field
$G(U)$ in (\ref{Eq2.8}). At the limit case, as $\varepsilon \to 0$,
the rate of change of these slow observables,
$\frac{d}{d\tau}\hat{v}_j$, is given by the right hand side of
(\ref{Eq5.1}), that is, by the average of $\nabla v_j(U) \cdot
G(U)$, the directional derivative of $v_j(U)$ along the vector field
$G(U)$, with respect of the Young invariant measure of the fast
dynamics (\ref{Eq3.1}). The need to average $\nabla v_j(U) \cdot
G(U)$  over the measure can be seen from the wild fluctuation in
\f{f3} of $\nabla v_3(U(\sigma)) \cdot G(U(\sigma))$, for the 3rd
Lax invariant $v_3$, along the trajectory $U(\sigma)$ of the fast
dynamics (\ref{Eq3.1}) that is shown in \f{f1}.
%
%
%
\begin{figure}[h]
\centerline{\psfig{figure=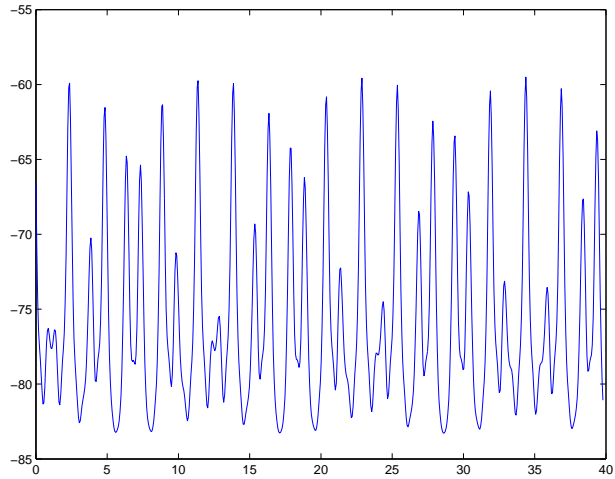,height=2.5in,height=2.5in}}
\caption{Evolution of $\nabla v_3(U(\sigma))\cdot G(U(\sigma))$
along the trajectory $U(\sigma)$ of  \f{f1}.}\label{f3}
\end{figure}

As we mentioned earlier in the introduction of section 6,  a typical
solution $U(\sigma)$ of the fast dynamics (\ref{Eq3.1}) was found to
be approximately a combination of two periodic motions, one fast and
the other  slower. The average value of $\nabla v_3(U(\sigma))\cdot
G(U(\sigma))$ along the trajectory $U(\sigma)$ of the fast dynamics
(\ref{Eq3.1}), from the time $\sigma = 0$ to the current time, is
shown in \f{f4}. Either the time over which the average is taken has
to be chosen carefully with respect to the approximate period of
$U(\sigma)$, or the average must be taken over a long time. In
\f{f4} the  circles indicate the averages at ``periodic" points (the
approximate fast period was estimated as 2.4868). Those dots
indicate that we can get a reasonable approximation for the Young
measure using a small number of approximate orbits.

\begin{figure}[h]
\centerline{\psfig{figure=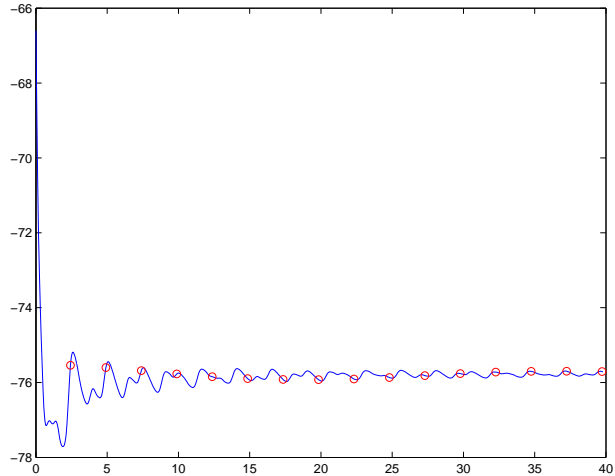,height=2.5in,height=2.5in}}
\caption{Averaged $\nabla v_3(U(\sigma))\cdot G(U(\sigma))$ along
the trajectory $U(\sigma)$ of \f{f1}.}\label{f4}
\end{figure}

We use this method of derivative evaluation of the slow observables,
by virtue of the average over the Young's measure in the right hand
side of (\ref{Eq5.1}), in order to integrate them forward in time by
simple forward Euler integration scheme with relatively large time
steps. Then we compare the performance, i.e. speedup and efficiency,
of this method to computing the slow observables through direct
numerical simulation of the full fast-slow system (\ref{nuequ}) with
$\nu = 10^{-4}$, where the latter will be considered as a ``true"
solution. Computing the slow observables via the Young's measure
approach  was done separately by first evaluating the time
derivatives of the slow observables, as it was described above, over
a fixed number of ``fast periods", and then advancing the slow
observables by taking an integration step in a forward Euler scheme
of either 3, 6, or 12 period length. We emphasize that the gain in
the efficiency of this method lies in the last step, i.e. in the use
of these relatively large time steps in advancing the slow
observables forward in time. The solutions for the second through
fourth slow observables, $v_2,v_3,v_4$, over 3,000 ``fast periods"
(a time interval of 7,460.4) are shown in \f{f5}. The slow
observables of the ``true" solution $U_\nu(t)$ (found by integrating
the full fast-slow system (\ref{nuequ}), with $\nu = 10^{-4}$, and
then computing the slow observables $v_j(U_\nu(t))$ directly from
it) are plotted as  solid
 lines, while the integration of the slow observables via the
Young's measure method, as described above, are plotted as  dotted
lines; but the difference (to plot accuracy) is {\bf not visible}.
\begin{figure}[h]
\centerline{\psfig{figure=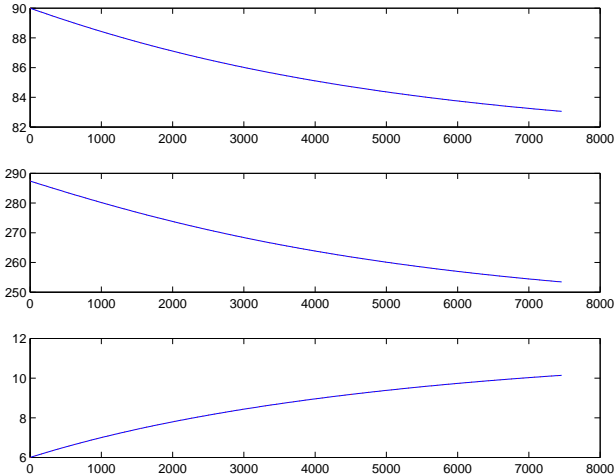,height=2.5in,height=2.5in}}
\caption{Behavior of the slow observables  $v_2,v_3$ and $v_4$ as
they are drifted by the slow diffusion.}\label{f5}
\end{figure}
We ran separate integrations  over either 1, 2, 3, or 4 fast
periods, in order to approximate the average in the right hand side
 of equation (\ref{Eq5.1}) with respect to the corresponding approximate
invariant Young measures of the fast system (\ref{Eq3.1}). These
separate integrations, or approximations of the Young measures, had
less than a 10\% effect on the solution errors. The dominant errors
were due to the Euler integration.  The errors in the slow
observables after 600 fast periods are shown in \t{T1}.
\begin{table}[h]
\center
\begin{tabular}{c|c c c }
  Euler step (periods) & $v_2$ & $v_3$ & $v_4$ \\ \hline
  3  & 0.45  & 0.47  &  0.68 \\
  6  &  1.3  &  2.1  &  1.8  \\
  12 &  3.9  & 10.8  &  4.1
\end{tabular}
\caption{Errors (scaled up by $10^3$) at $t = 1,492.08$ (600 fast
periods) in the slow observables  integrated using approximations to
Young measure averages.} \label{T1}
\end{table}
The speedup and efficiency  of the Young measure integration scheme,
in comparison to direct numerical simulation of the full fast-slow
system (\ref{nuequ}), depends on the value of $\nu$. The cost of
evaluating the time derivative of the slow observables via the Young
measure averaging in the right hand side of (\ref{Eq5.1}) is
determined solely by the accuracy needed. Moreover, the number of
integration steps needed for advancing the slow observables forward
in time  is also determined solely by the accuracy needed. However,
if we were to integrate equation (\ref{nuequ}) directly, the number
of time steps needed would be proportional to the integration
interval which is proportional to $1/\nu$.  If, for example, we
needed 50 integration steps for one approximate period, $\nu =
0.001$, and we wanted to track the dissipative (diffusion) terms
until they had reduced by a factor of $10^{-5}$ we would need to
integrate from $t = 0$ to $t = 11,500$ (approximately).  Since each
approximate period is 2.4868 we would need about 4,600 periods, or
23,000 integration steps for a direct numerical integration of
(\ref{nuequ}).  If we used a higher-order integration method (rather
than the first order Euler method used in this illustration) at the
projective level, namely for advancing the slow observables in time
after evaluating their time derivatives, we could probably integrate
the decaying slow observables  with much less than 100 projective
steps - each of which requires one integration around of the fast
system using 50 steps, for a total of 5,000 steps, or a speedup of
better than 4 to 1. If $\nu$ were reduced by one tenth, the speedup
would increase to better than 40 to one.

As mentioned, while computing the invariant measures of the fast
system (\ref{Eq3.1}) (or equivalently of (\ref{nuequ}) with $\nu=0$)
we actually simulate the entire evolution of the limit system of
(\ref{nuequ}), as $\nu \to 0$, via the supports of these measures -
tori in our example. The decay of the tori of the full fast-slow
system (\ref{nuequ}), with $\nu >0$, can be seen in \f{fm} which
shows the trajectories after 600, 1200, 1800, 2,400, and 3,000 fast
periods of integration.

\begin{figure}[h]
\centerline{\psfig{figure=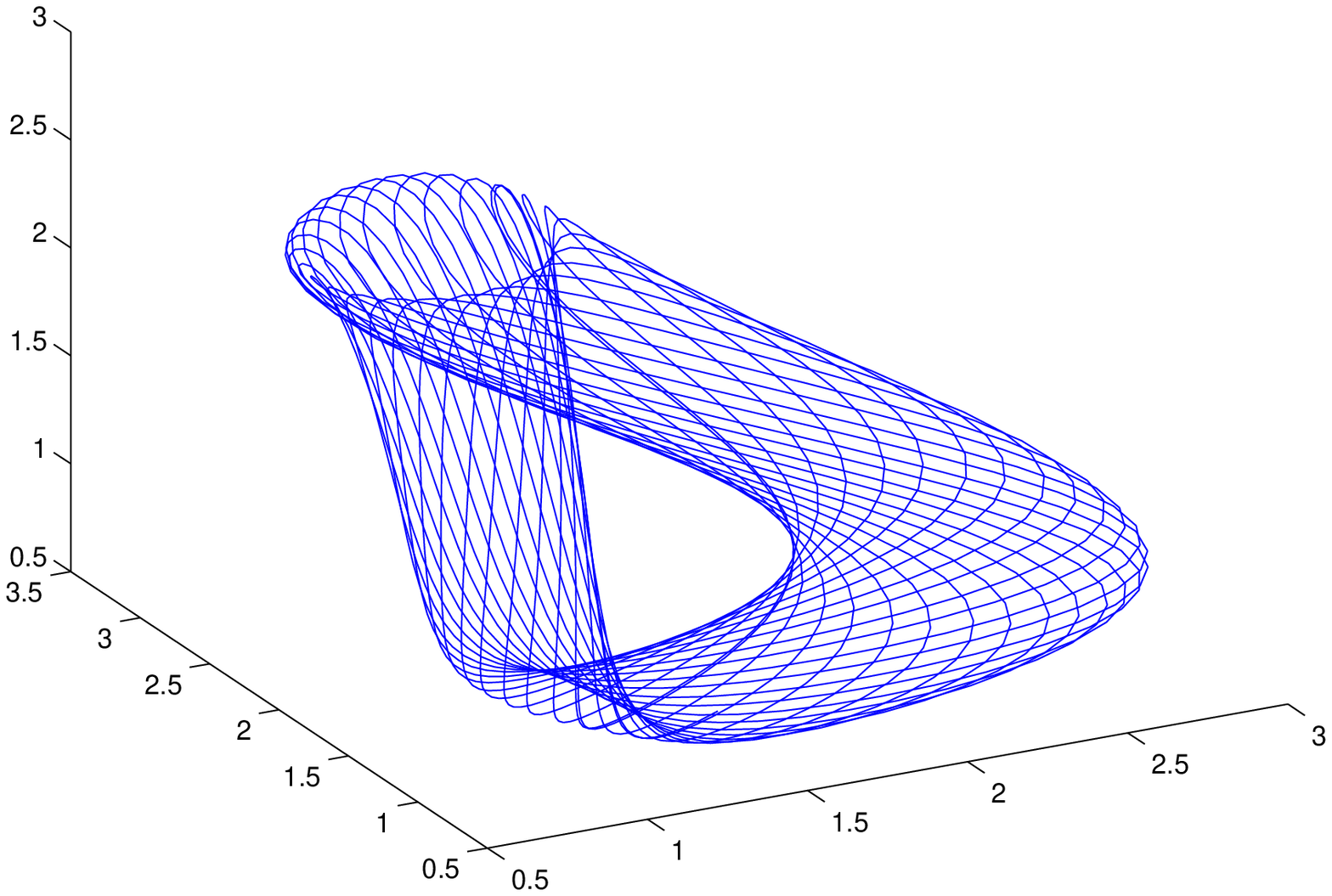,height=1.5in,height=1.5in}
 \psfig{figure=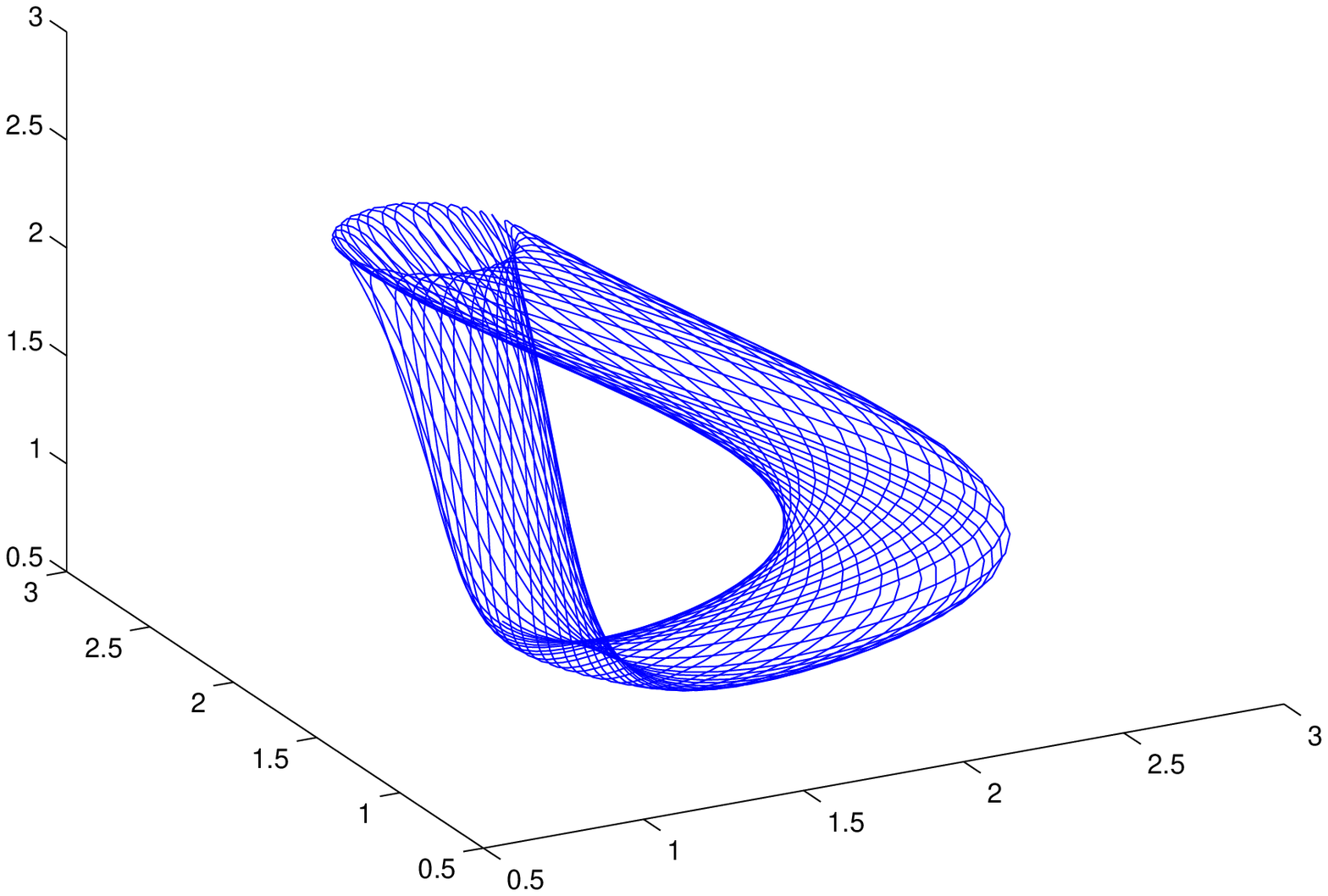,height=1.5in,height=1.5in}
\psfig{figure=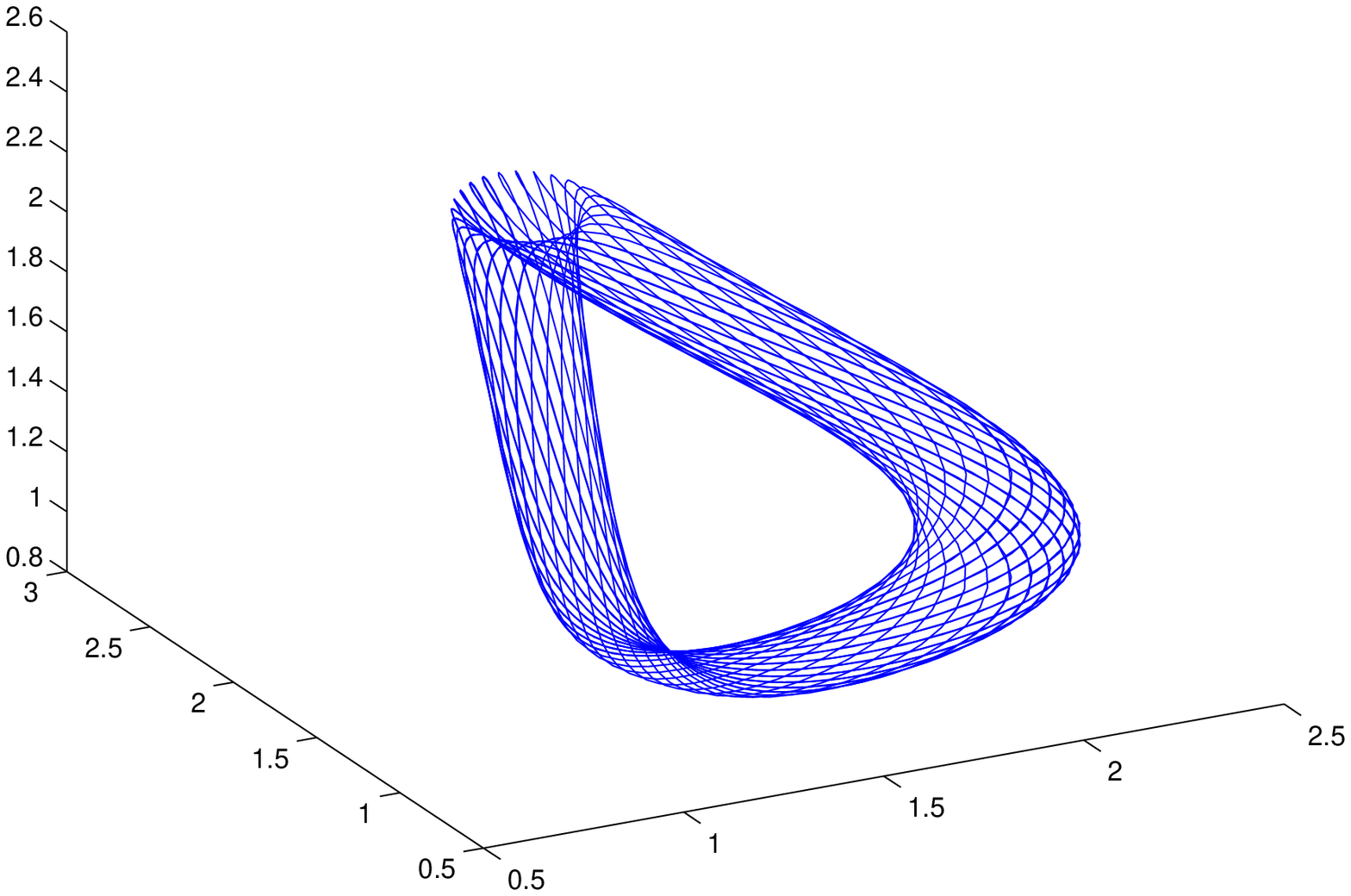,height=1.5in,height=1.5in}}
\centerline{\psfig{figure=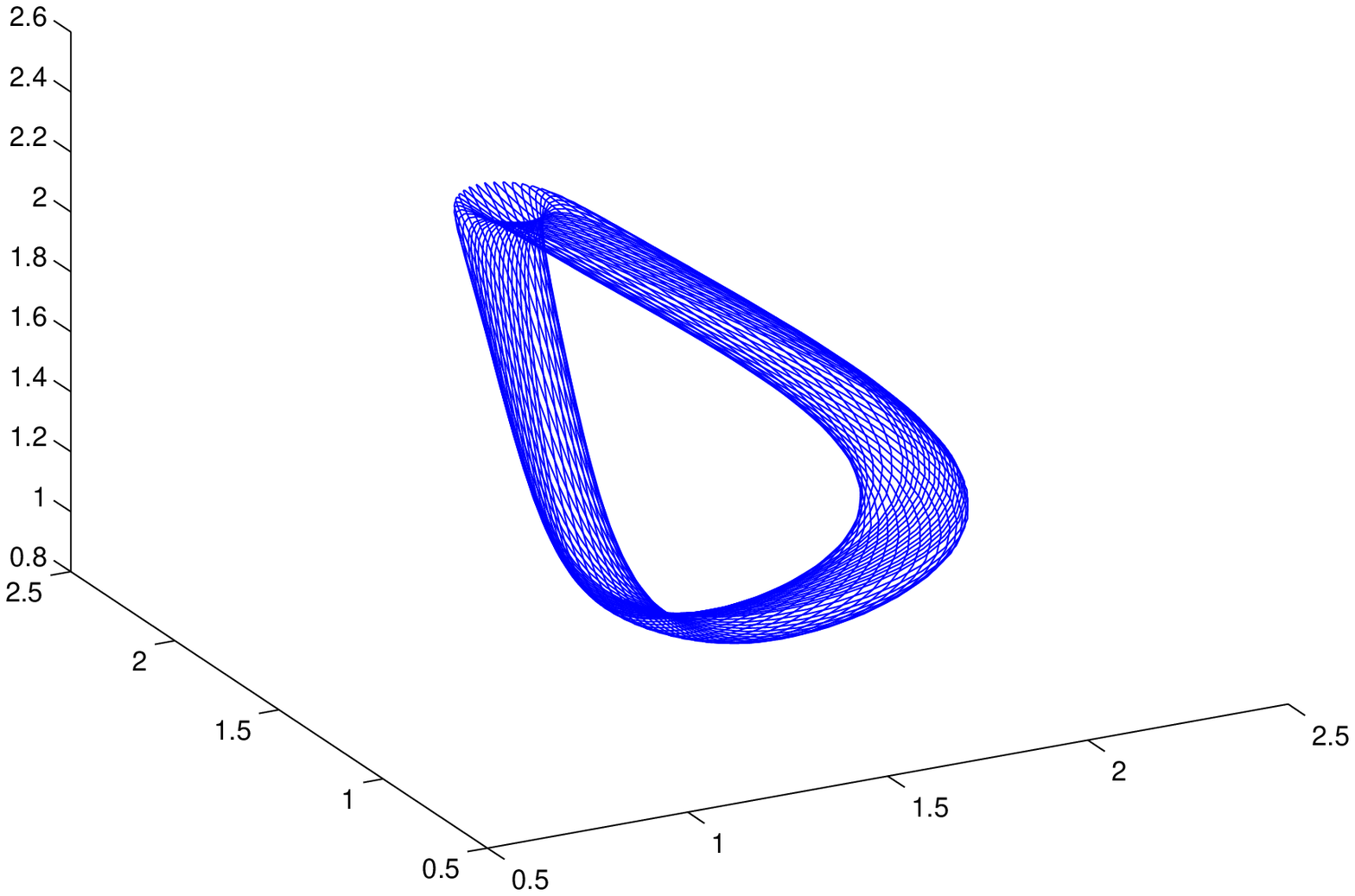,height=1.5in,height=1.5in}
\psfig{figure=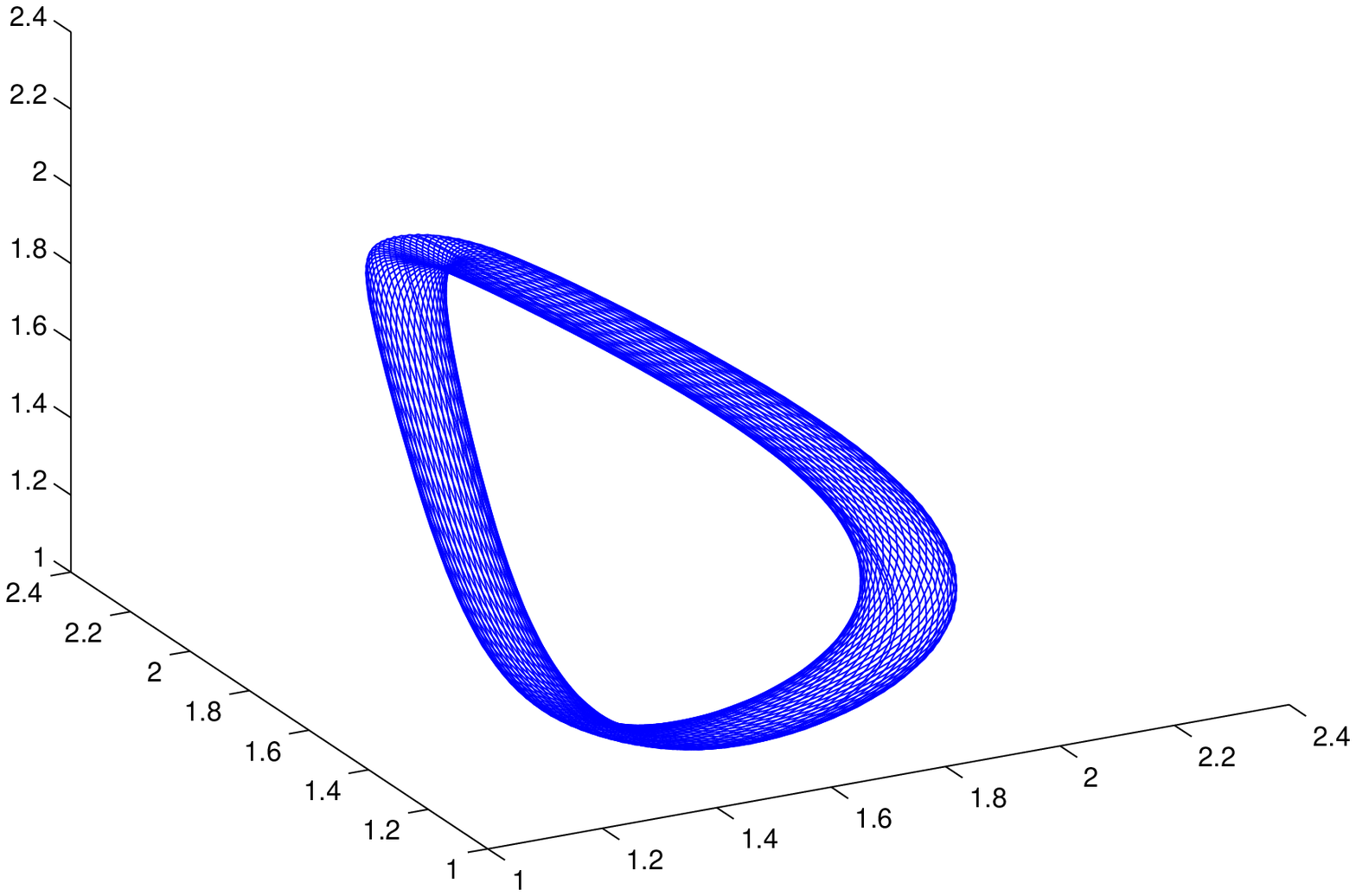,height=1.5in,height=1.5in}
\psfig{figure=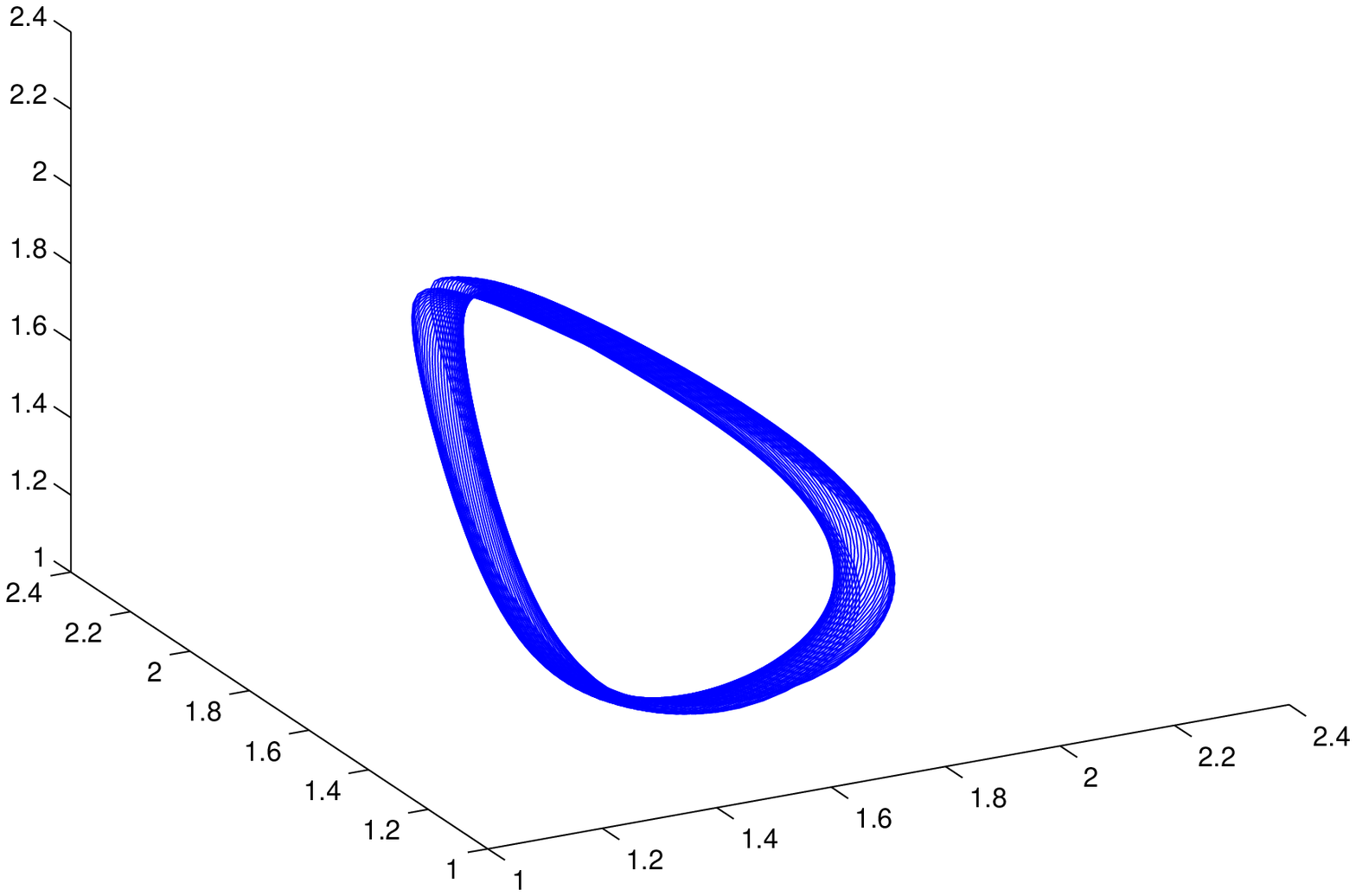,height=1.5in,height=1.5in}}
\caption{Evolution of Tori of system (\ref{nuequ}) for initial
condition [1 1 1 1 4 1].}\label{fm}
\end{figure}


\subsection{The second approach - the Equation-Free method}

The equation-free method \cite{KGHKRT} utilizes the explicit
knowledge of the numerical solution $U_\nu (t)$ of (\ref{nuequ}),
and consequently of the observables $v_j(U_\nu(t))$,
$j=1,2,\cdots,\frac{N}{2}+1$, for short interval of time in order to
approximate, here by numerical differencing, the time derivative of
 $v_j(U_\nu(t))$ in order to project forward their values. The
 method is called equation-free since it does not use explicitly any
equation of motion for the slow observables. However, it assumes
that the evolution of $v_j(U_\nu(t))$ is slow, and that there is an
underlying process, or equation, that governs their evolution, which
we do not use explicitly. (A valuable research direction would be to
develop a method that identifies, possibly numerically, the slow
observables in case they are not prescribed.)

 For small values of $\nu > 0$ the evolution of some of the slow
observables, $v_j(U_\nu(t))$, $j=1,2,\cdots,\frac{N}{2}+1$, along
solutions $U_\nu(t)$ of (\ref{nuequ}) can exhibit a slight
oscillatory behavior. For example, \f{f6} shows the early behavior
of $v_3(U_\nu(t))$ along the solution $U_\nu(t)$ of equation
(\ref{nuequ})  for the initial condition [1 1 1 1 4 1]. It exhibits
a fine-scale oscillation. However, the values at the ``periodic"
points are reasonably smooth so we can use those points in order to
evaluate the time derivative of the slow observables,
$v_j(U_\nu(t))$, by numerical differencing. Then we use the
numerical time derivatives of $v_j(U_\nu(t))$   to project them
forward in time. The results are illustrated in \f{f7}. We show the
results of projective integration of the slow observables (dotted
lines) versus the evolution of the slow observables along a true
solution of (\ref{nuequ}) (solid line) for the three slow
observables $v_2,v_3$ and $v_4$. Again, and as in the Young measure
method, the results are indistinguishable to plot accuracy.

\begin{figure}[h]
\centerline{\psfig{figure=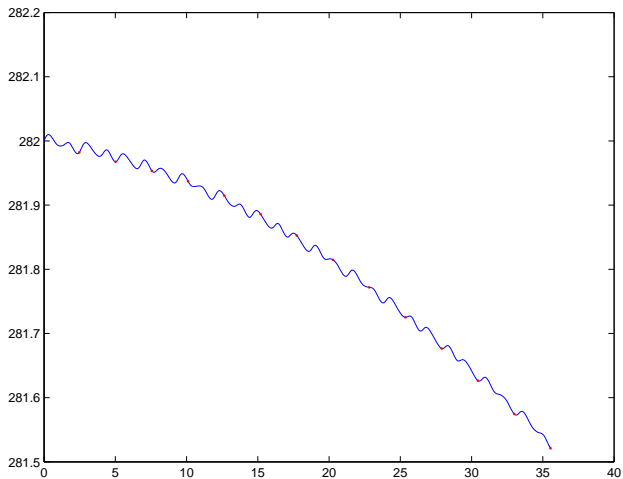,height=2.5in,height=2.5in}}
\caption{Slow observable $v_3(U_\nu(t))$ for initial condition [1 1
1 1 4 1].}\label{f6}
\end{figure}

\begin{figure}[h]
\centerline{\psfig{figure=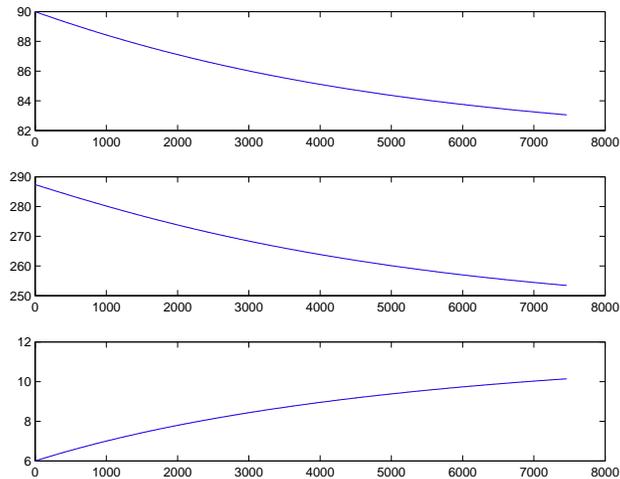,height=2.5in,height=2.5in}}
\caption{Behavior of the slow observables $v_2,v_3$ and $v_4$ as
they are drifted by the slow diffusion using the Equation-free
integration.}\label{f7}
\end{figure}

\t{T2} shows the errors in the slow observables at $t = 1,492.08$
using projective Euler integration by evaluating the time derivative
from the slope of the chord over one period.

\begin{table}[h]
\center
\begin{tabular}{c|c c c }
  Euler step (periods) & $v_2$ & $v_3$ & $v_4$ \\ \hline
  3  & 0.8  & 3.6  & 0.5  \\
  6  & 1.9  & 5.9  & 1.8  \\
  12 & 4.2  & 14.2  & 3.9
\end{tabular}
\caption{Errors (scaled up by $10^3$) at $t = 1,492.08$ (600 fast
periods) in the slow observables  integrated using the Equation-free
method over one period.}\label{T2}
\end{table}

\section*{Acknowledgements} \noindent
The research of Z.A., I.G.K., M.S. and E.S.T. was supported  by a
grant from the United States - Israel Binational Science Foundation
(BSF).

\end{document}